\documentclass[12pt]{amsart}
\usepackage{amssymb,latexsym, amsthm,epsf}%, makeidx}
\newtheorem{thm}{Theorem}[section] %the resolution could also be [subsection]
\newtheorem{defn}[thm]{Definition}

\newtheorem{const}[thm]{Construction}

\newtheorem{lem}[thm]{Lemma}

\newtheorem{prop}[thm]{Proposition}
\newtheorem{rem}[thm]{Remark}

\def\N{{\mathbb N}}
\def\R{{\mathbb R}}

\def\({\left(}
\def\){\right)}
\long\def\forget#1\forgotten{}

\newif \iffurther %Give \section{Further ideas}
\furtherfalse

\newif \iffurther %Give \section{Further ideas}
\furtherfalse

\newif\ifXY % turns XY version on/off
\XYfalse    % Turn it off

\ifXY
\usepackage{xy}
\fi
\ifXY
\xyoption{all}
\fi

\begin{document}

\title[Chromatic properties of point configurations]{Chromatic properties of generic planar configurations of points}

\author[Roland Bacher and David Garber]{Roland Bacher$^1$ and David Garber$^2$}

\stepcounter{footnote} 
\footnotetext{Support from the Swiss National Science Foundation is 
gratefully acknowledged.}
\stepcounter{footnote} 
\footnotetext{Partially supported by the Chateaubriand postdoctoral
  fellowship funded by the French government.} 

\address{Institut Fourier, BP 74, 38402 Saint-Martin D'Heres CEDEX,
  France}
\email{\{bacher,garber\}@mozart.ujf-grenoble.fr}
\date{\today}

\begin{abstract} 
We study the Orchard relation defined in \cite{bacher} for
generic configurations of points in the plane (also called order
types). 
We introduce infinitesimally-close points 
and analyse the relation of this notion with the Orchard
relation. 

The second part of the paper deals with
monochromatic configurations (for the Orchard relation). We
give the complete list of all monochromatic configurations 
up to $7$ points
and present some constructions and families of monochromatic
configurations.
\end{abstract}

\maketitle

\section{Introduction}
A finite set $\mathcal P = \{P_1,\cdots, P_n\}$ of $n$ points 
in the oriented affine plane $\R^2$ is a {\it generic configuration} 
if three points in $\mathcal P$ are never collinear. 
Two generic configurations of $n$ points $\mathcal P^1$ and 
$\mathcal P^2$ are {\it isotopic} if they can be joined by a continuous
curve of generic configurations. 

Two generic configurations ${\mathcal P}^1$ and ${\mathcal P}^2$ are
{\it isomorphic} if there exists a bijection 
$\varphi: {\mathcal P}^1\longrightarrow{\mathcal P}^2$ such that 
the two triangles having vertices 
$P,Q,R\in{\mathcal P}^1$ and $\varphi(P),\varphi(Q),\varphi(R)\in{\mathcal 
P}^2$ induce either always identical or always opposite orientations
for all triplets
$\{ P,Q,R \}$ of points in ${\mathcal P}^1$. In the former case we call 
${\mathcal P}^1$ and ${\mathcal P}^2$ {\it orientedly isomorphic}.
The (oriented) isomorphism classes of all generic configurations are
called {\it (oriented) order types} by some authors, especially in
Computer Sciences (see \cite{aak} and \cite{ak}).

Isotopic configurations are of course orientedly isomorphic. 
We ignore to what extend the converse holds.

A line $L \subset \R^2$ {\it separates} two points $P,Q \in \R^2\setminus L$ if
$P$ and $Q$ are in different connected components of $\R^2 \setminus L$. 
Given a generic configuration $\mathcal P$, we denote by $n(P,Q)$
the number of separating lines defined by pairs of points in $\mathcal P 
\setminus \{P,Q\}$.

\begin{defn}[Orchard relation] We set
$P \sim Q$ if we have $n(P,Q) \equiv (n-3) \pmod 2$ for two distinct points 
$P,Q\in{\mathcal P}$ of a generic configuration ${\mathcal P}$.
\end{defn}

The following result is then a special case of a more general
fact, see \cite{bacher}.

\begin{thm}[Orchard Theorem]\label{orchard_thm}
The Orchard relation is an equivalence relation on the set of points
$\mathcal P$
of a generic configuration, which consists of at most two classes.
\end{thm}

The name of this result comes from the fact that it yields a
canonical rule to
plant trees of two species at prescribed generic locations in an
orchard.

We call the induced (generally non-trivial)
partition the {\it Orchard partition} of $\mathcal P$.
Orchard partitions yield invariants for studying generic
configurations of points in the affine plane. Iterative use of the Orchard
partition, i.e. dividing the set of points into
Orchard classes, and iterating this on each class,
produces an invariant which distinguishes many pairs of 
non-isomorphic generic configurations and which is easy to compute
and handle.
  
\medskip

The present paper is devoted to the illustration
of the planar Orchard Theorem applied to
generic configurations of $n$ points in the affine plane.
Such configurations arise 
naturally and many features of them have been considered
by other authors. This paper yields another such contribution
devoted to  Orchard properties.

Configurations with trivial Orchard partition (monochromatic
configurations) are especially noteworthy to study since they 
are the \lq\lq atoms'' of the theory.
We will describe a few infinite monochromatic families 
and constructions involving monochromatic configurations.

We introduce and study also some properties of a notion
which we call {\it infinitesimal-closedness}. This property is 
useful for understanding some monochromatic families and has
also some independent interest.

\medskip

The paper is organized as follows. Section \ref{proofs} recalls the proof of
the Orchard Theorem in the planar case. We also recall and prove 
the corresponding flip proposition. 
Section \ref{examples} lists all generic configurations 
having up to $6$ points, together with their Orchard partitions. 
Section \ref{infi_close_section} deals with infinitesimally-close
points and the effect on the Orchard relation after deleting 
two such points.

The remaining  part of the paper is centered on
monochromatic configurations. Section \ref{mono_config}
gives the complete list of monochromatic configurations
up to $7$ points. Section \ref{mono_construct} deals with some
constructions preserving monochromatic configurations.
Section \ref{mono_families} describes some monochromatic
families. 

Section \ref{statistics} uses a data-base of Aichholzer,
Aurenhammer and Krasser for collecting some statistics related to the
Orchard partitions on the set of all
generic configurations having at most $9$ points.

\section{Proofs of the Orchard Theorem and of the flip property for
  planar configurations}\label{proofs}

In this section we prove the planar version of the
Orchard Theorem and the
flip property along the lines of \cite{bacher}. The planar
setting allows a slight simplification.

\begin{figure}[h]
\epsfysize=5cm
\centerline{\epsfbox{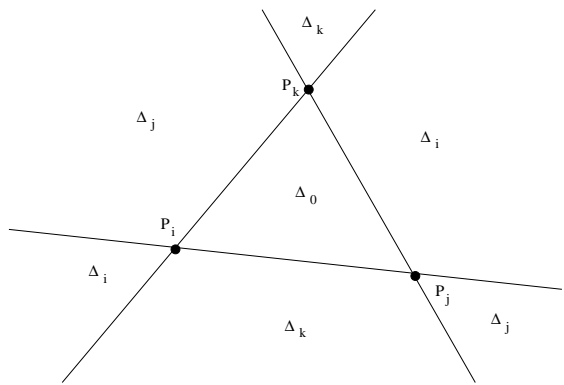}}
\caption{Division of the projective
plane into four triangles $\Delta_0 \cup \Delta_i \cup \Delta_j \cup \Delta_k$}\label{orchard_proof}
\end{figure}

Three points $P_i,P_j,P_k$ of a generic planar configuration 
$\mathcal P = \{P_1 , \cdots P_n \}$ define three lines which 
subdivide the projective plane into four
triangles: $\Delta_0,\Delta_i,\Delta_j,\Delta_k$, as illustrated in
Figure \ref{orchard_proof}.
  
Denote by $\alpha_i$ the number of lines connecting two
points in $\mathcal P \setminus \{ P_i,P_j,P_k \}$ and separating (in
the affine plane of course)
$P_i$ from both points $P_j$ and $P_k$. The numbers
$\alpha_j$ and $\alpha_k$ are defined analogously. 

For $* \in \{ 0,i,j,k \}$ denote by $\sigma_*$ the number of points in 
$\mathcal P \setminus \{ P_i,P_j,P_k \}$ which are in the interior
of the triangle $\Delta_*$.

\begin{lem}\label{orchard_lemma} We have
$$n(P_i,P_j) = \alpha_i+\alpha_j + \sigma_0 + \sigma_k$$
$$n(P_j,P_k) = \alpha_j+\alpha_k + \sigma_0 + \sigma_i$$
$$n(P_i,P_k) = \alpha_i+\alpha_k + \sigma_0 + \sigma_j$$
\end{lem}

\begin{proof}
We prove the first formula. The others follow by symmetry.

For a line $L$ separating $P_i$ from $P_j$ we have two cases: 
$P_k \not\in L$ or $P_k \in L$.

In the first case, $P_k \not\in L$, the line $L$ separates either 
$P_i$ from $P_k$ or $P_j$ from $P_k$ (see e.g. line $L_1$ 
in Figure \ref{lemma_fig}). 
Such a line $L$ yields hence a contribution of $1$ to either $\alpha_i$
or $\alpha_j$. 

If $P_k \in L$, then $L$ is defined by the point $P_k$ and by another
point $Q$ which belongs either to $\Delta_0$ or to
$\Delta_k$ (consider line $L_2$ or $L_3$ in Figure \ref{lemma_fig}).

\begin{figure}[!h]
\epsfysize=7cm
\centerline{\epsfbox{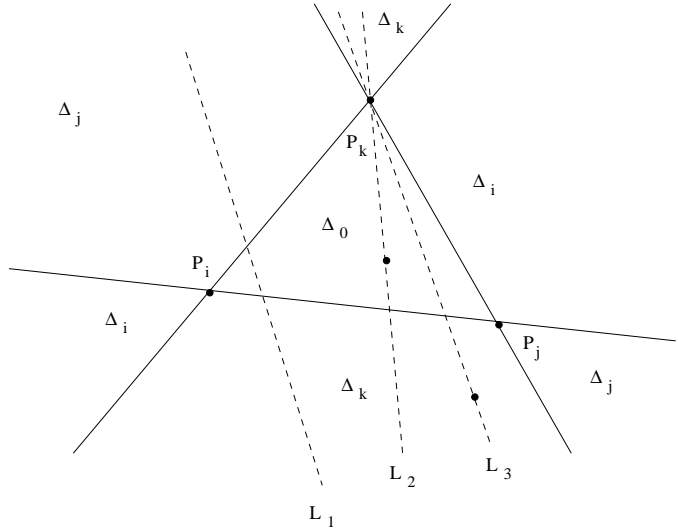}}
\caption{A few lines in a generic configuration}\label{lemma_fig}
\end{figure}                                                                 
 
\end{proof}

\begin{proof}[Proof of Theorem \ref{orchard_thm}]
Reflexivity and symmetry of the Orchard relation are obvious.
We have to establish transitivity.

Let $P_i,P_j,P_k$ be three points such that $P_i \sim P_j$ and   
$P_j \sim P_k$. We have to show that $P_i \sim P_k$.

Since $P_i \sim P_j$ and $P_j \sim P_k$ we have
$n(P_i,P_j) \equiv (n-3) \pmod 2$ and 
$n(P_j,P_k) \equiv (n-3) \pmod 2$. 
Adding the equations 
$$n(P_i,P_j) = \alpha_i+\alpha_j + \sigma_0 + \sigma_k \equiv (n-3) \pmod 2$$
$$n(P_j,P_k) = \alpha_j+\alpha_k + \sigma_0 + \sigma_i \equiv (n-3) \pmod 2$$
obtained by Lemma \ref{orchard_lemma}, we get
$$\alpha_i+\alpha_k+\sigma_i+\sigma_k \equiv 0 \pmod 2\ .$$
The interiors of the four triangles $\Delta_0,\Delta_i,\Delta_j$ and
$\Delta_k$ contain all points of 
${\mathcal P}\setminus\{
P_i,P_j,P_k\}$ by genericity and we have hence
$$\sigma_i+\sigma_j+\sigma_k+\sigma_0 = n-3\ .$$
We get therefore
$$n(P_i,P_k) = \alpha_i+\alpha_k + \sigma_0 + \sigma_j \equiv
\alpha_i+\alpha_k + \sigma_i + \sigma_k+(n-3)\equiv (n-3) \pmod 2$$
which shows that $P_i \sim P_k$.  

The fact that the Orchard relation has at most two classes 
follows from the implication $P_i\not\sim P_j$ and $P_j\not\sim P_k
\Longrightarrow P_i\sim P_k$ (the proof of which is, up to a minor
change, as above).
\end{proof}

\medskip

A {\it flip} is the most elementary move relating non-isomorphic 
generic configurations:
\begin{defn}[Flip]
Two generic configurations $\mathcal P ^1, \mathcal P ^2$ of $n$ points
are related by a flip 
if there exists a continuous path of configurations 
$\mathcal P(t)$, $-1 \leq t \leq 1$, such that:
\begin{enumerate}
\item $\mathcal P (-1) = \mathcal P ^1$ and $\mathcal P (1) = \mathcal P ^2$.
\item $\mathcal P (t)$ is generic for all $-1 \leq t \leq 1$ except 
for $t=0$. The configuration $\mathcal P (0)$ has exactly three
aligned points, each crossing transversally at $t=0$ 
the line spanned by the two other points.
\end{enumerate}
\end{defn}

\begin{figure}[!h]
\epsfysize=4cm
\centerline{\epsfbox{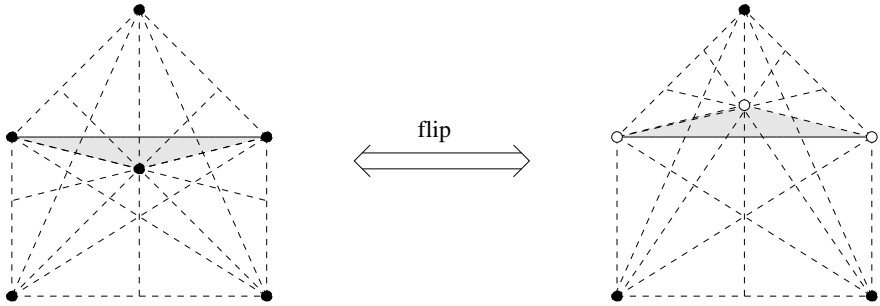}}
\caption{Example of a flip}\label{flip_exam}
\end{figure}

Geometrically, a flip transforms an almost \lq\lq flat'' 
triangle formed by three points
of $\mathcal P$ which are nearly aligned 
into its \lq\lq mirror'', see Figure \ref{flip_exam} for an illustration. We
leave it to the reader to show that every pair of generic
configurations with $n$ points
can be related by a path involving only isotopies and a finite number
of flips. 

A flip affects the Orchard
relation only locally as shown by the following proposition:

\begin{prop}[Flip Proposition]\label{flip_prop}
Let $\mathcal P=\mathcal P (-1)= \{P_1, \cdots, P_n\}$ and 
$\mathcal P '=\mathcal P (1)= \{P'_1, \cdots, P'_n\}$ be two generic configurations
related by a flip involving three points $R,S,T\in \mathcal P$ 
corresponding to 
$R',S',T'\in \mathcal P '$.
For any two points $P,Q$ we have:
\begin{enumerate} 
\item if $P,Q \in \{ R,S,T \}$ or 
$P,Q \in \mathcal P \setminus \{ R,S,T \}$, then:
$$P \sim Q \quad {\rm if\ and\ only\ if} \quad P' \sim' Q'$$
\item if $P \in \{ R,S,T \}$ and $Q \in \mathcal P \setminus \{ R,S,T \}$
then:
$$P \sim Q \quad {\rm if\ and\ only\ if} \quad P' \not\sim' Q'$$
\end{enumerate} 
\end{prop}

\begin{proof} 
A point $V(t)\in \mathcal P (t)\setminus \{R(t),S(t),T(t)\}$ is never
aligned with two other points of $\mathcal P (t)$. This shows the equality
$$n(P(1),Q(1))=n(P(-1),Q(-1))$$ 
for $P(t),Q(t)\in \mathcal P (t)\setminus \{R(t),S(t),T(t)\}$ or 
$P(t),Q(t)\in \{R(t),S(t),T(t)\}$.

On the other hand, $V(t)\in \{R(t),S(t),T(t)\}$ will be involved
transversally in the unique alignement $\{R(0),S(0),T(0)\}$ 
during the flip. For such a point $V(t)$ and $P(t)\in
\mathcal P (t)\setminus \{R(t),S(t),T(t)\}$ we have hence
$$n(V(1),P(1))=n(V(-1),P(-1))\pm 1$$
which proves the result.
\end{proof}

\section{All generic configurations up to $6$ points}\label{examples}

This section contains the complete list (up to unoriented
isomorphisms) of all generic configurations having at most $6$
points together with the corresponding Orchard partitions.  

We represent the two classes by black and white vertices. The
choice of the black class is of course arbitrary: In case of classes
having different cardinalities, black is used for
the class containg more elements. In case of a draw, 
the choice of the black class is arbitrary. 

\subsection{Generic configurations with up to $5$ points}

\begin{figure}[!h]
\epsfysize=6cm
\centerline{\epsfbox{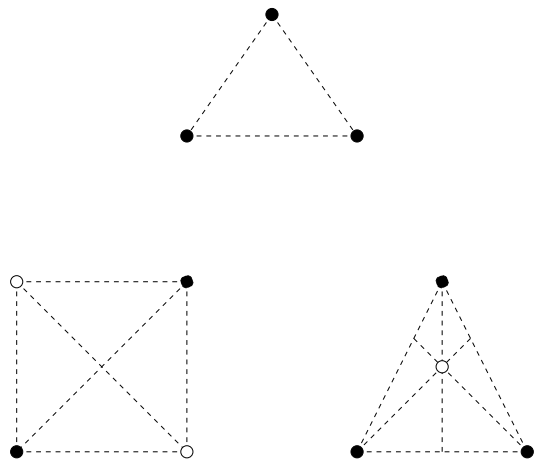}}
\caption{All generic configurations with $3$ or $4$ points}\label{part_34points}
\end{figure}

\begin{figure}[!h]
\epsfysize=7cm
\centerline{\epsfbox{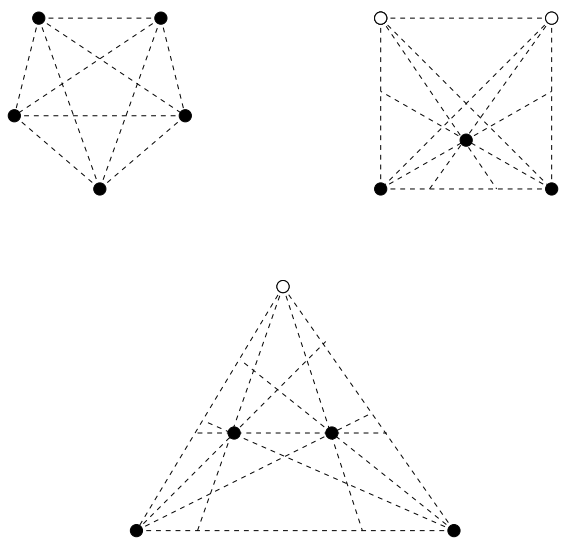}}
\caption{All generic configurations with $5$ points}\label{part_5points}
\end{figure}

Figures \ref{part_34points} and \ref{part_5points}
present all generic configurations (up to 
non-oriented isomorphism) having at most $5$ points,
together with their Orchard partitions. Since all these classes
admit an orientation-reversing symmetry, these classes coincide with
oriented isomorphism classes.

\newpage

\subsection{Generic configurations of $6$ points}

\begin{figure}[h]
\epsfysize=7cm
\centerline{\epsfbox{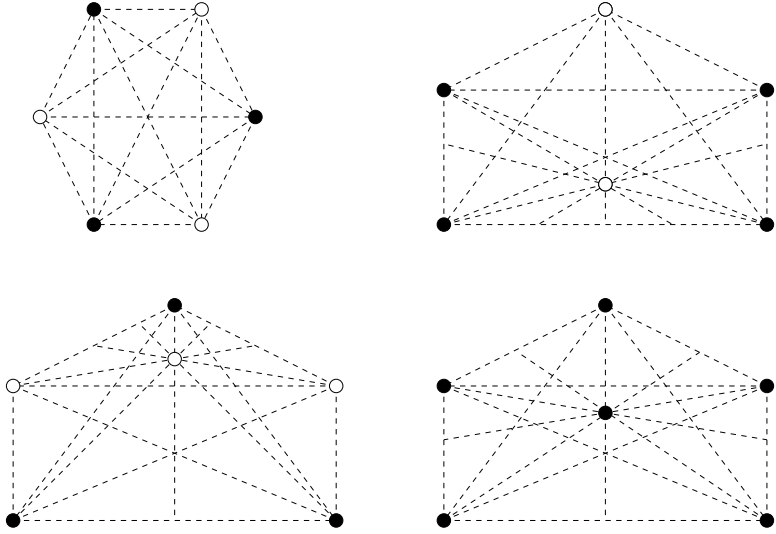}}
\caption{Generic configurations of $6$ points with a convex hull of size
  $5$ or $6$}\label{part_6points_1}
\end{figure}

\begin{figure}[h]
\epsfysize=15cm
\centerline{\epsfbox{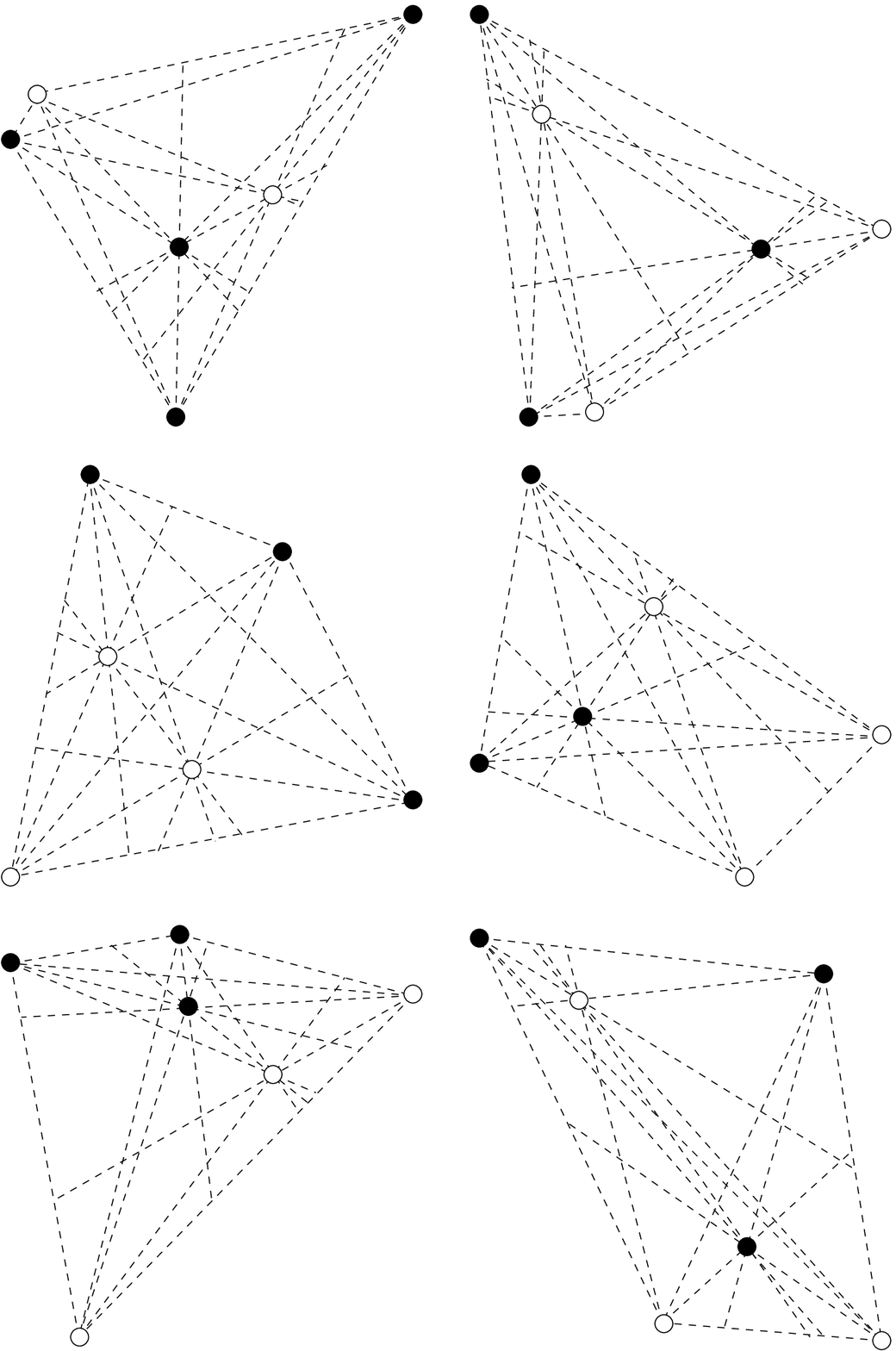}}
\caption{Generic configurations of $6$ points with a convex hull of
  size $4$}\label{part_6points_2}
\end{figure}

\begin{figure}[h]
\epsfysize=16cm
\centerline{\epsfbox{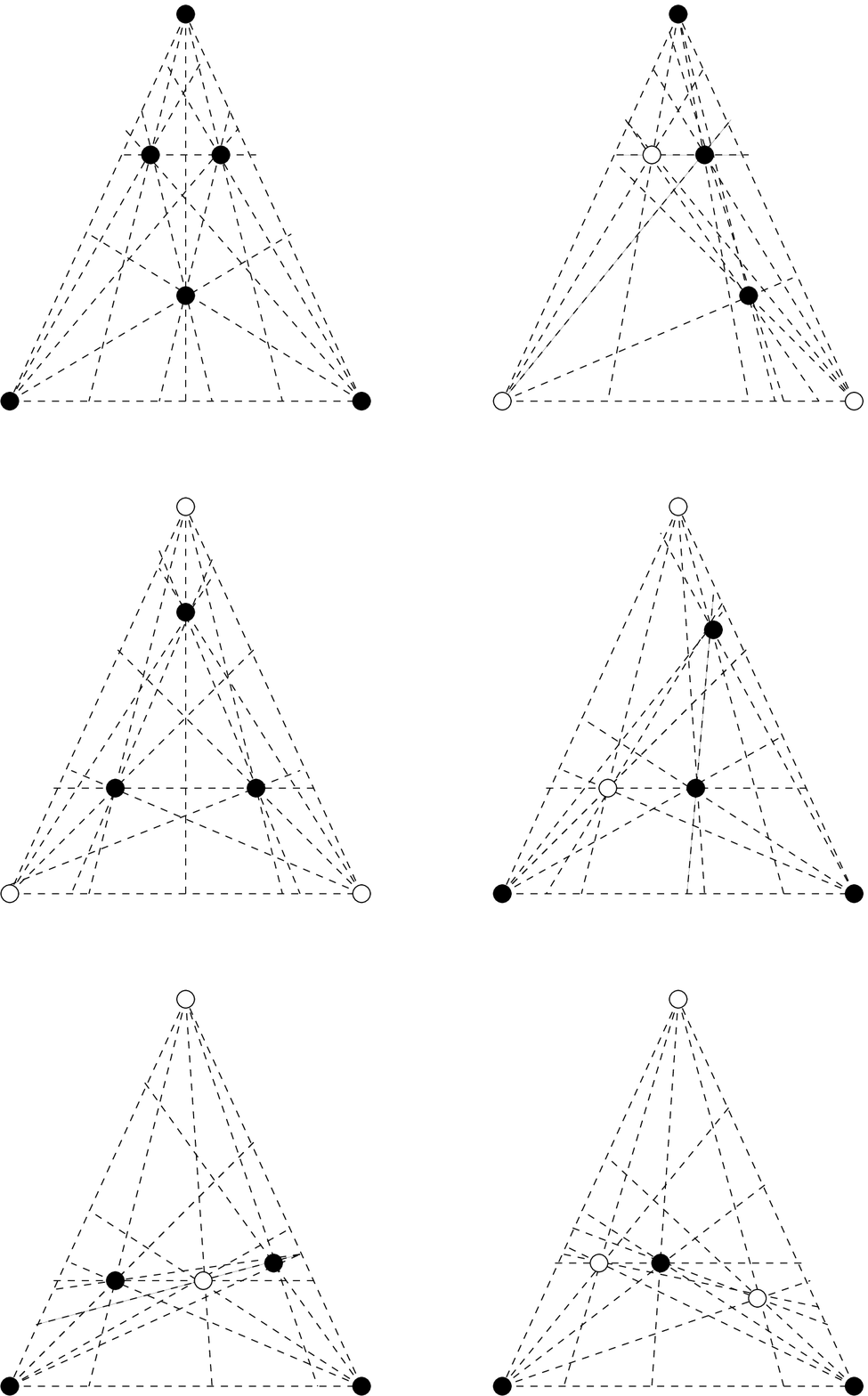}}
\caption{Generic configurations of $6$ points with a convex hull of size $3$}\label{part_6points_3}
\end{figure}

Figures \ref{part_6points_1}, \ref{part_6points_2} and
\ref{part_6points_3} contain all non-oriented 
isomorphism classes of generic configurations with $6$
points (sorted by the size of their convex hull), together
with their Orchard partitions. 

\newpage

\newpage

\section{Infinitesimally-close points}\label{infi_close_section}
In this section, we study configurations containing a pair 
$P,Q$ of points with $n(P,Q)=0$ (no separating lines). 
Such pairs of points will be useful for
constructing monochromatic configurations.

\begin{defn}
Two points $P,Q$ are {\it infinitesimally-close} if $n(P,Q)=0$. 
\end{defn}

\begin{rem}\label{infi_close_rem} Infinitesimally-close points are
Orchard-equivalent if the total number of points in the
configuration is odd and they are not equivalent
otherwise.
\end{rem}

The following proposition describes the change in the Orchard relation
induced by deleting two infinitesimally-close points:

\begin{prop}\label{infi_close_prop}
Let
$R,S \in \mathcal P$ be two infinitesimally-close points of a generic
configuration $\mathcal P$ having $n$ points. Denote by $L$ the
line containing $R$ and $S$ and by 
$\mathcal P' = \mathcal P \setminus \{ R,S \}$
the subconfiguration of $\mathcal P$ after deletion of $R$ and $S$. 
\begin{enumerate}
\item If $P,Q \in \mathcal P \setminus \{ R,S \}$ are not separated by $L$, then
$$P \sim_{\mathcal P} Q  \quad {\rm if\ and\ only\ if} \quad 
P \sim_{\mathcal P'} Q\ .$$
\item If $P,Q \in \mathcal P \setminus \{ R,S \}$ are separated by $L$, then
$$P \sim_{\mathcal P} Q \quad {\rm if\ and\ only\ if} 
\quad P \not\sim_{\mathcal P'} Q\ .$$
\end{enumerate}
\end{prop}

\begin{proof}
If $P,Q$ are not separated by $L$, 
then $n_{\mathcal P}(P,Q)$ and $n_{\mathcal P'}(P,Q)$ have the
same parity. Indeed, if a separating line between $P$ and $Q$ is 
generated by $R$ (say) and a second point 
$T \in {\mathcal P} \setminus \{ R,S,P,Q \}$, 
then the line generated by $S$ and $T$ separates $P$ and
$Q$ too. The number of such lines is hence even.

If the line $L$ spanned by $R$ and $S$ separates two points
$P,Q\in{\mathcal P}\setminus\{R,S\}$ then the number
$n_{\mathcal P}(P,Q)-n_{\mathcal P'}(P,Q)$
is odd since it contains, together with all pairs of lines
mentioned above, also the line $L$.
\end{proof}

\begin{rem}
Generic configurations consisting of few 
points have generally pairs of infinitesimally-close points.
Generic configurations having a huge number of points however  
have only rarely infinitesimally-close points.
\end{rem}

\subsection{Infinitesimal contractions and infinitesimally
equivalent configurations}

Let $\mathcal P$ be a configuration with two infinitesimally-close
points $R$ and $S$. This
implies that the two configurations ${\mathcal P}\setminus\{R\}$ 
and ${\mathcal P}\setminus\{S\}$,
obtained by deleting either $R$ or $S$ from ${\mathcal P}$,
are isotopic. We call the isotopy class $\mathcal P '$ 
of $\mathcal P \setminus\{R\}$
the {\it infinitesimal contraction} of $\mathcal P$ along $[R,S]$.
Two configurations $\mathcal P,\ \mathcal P '$ related through 
an infinitesimal contraction are {\it infinitesimally related}.
Two configurations $\mathcal P ^1$ and $\mathcal P ^2$ which can be
joined by a path of infinitesimally related configurations are
{\it infinitesimally equivalent} and the set of all 
configurations infinitesimally equivalent to $\mathcal P$ forms of
course an equivalence class. Nothing interesting can be said 
concerning Orchard properties of such a class since removal of a 
point changes the Orchard relation dramatically.
It has however a nice property given by the following result.

\begin{thm} 
Each equivalence class of infinitesimally related
configurations contains a unique (up to isomorphism) minimal
representative. The minimal representative of a configuration $\mathcal
P$ can always be constructed as a subconfiguration
reached through a succession of infinitesimal
contractions starting at $\mathcal P$.
\end{thm}

\begin{proof} 
We have to prove unicity of a minimal representative.
Choose a minimal sequence $\gamma(0),\gamma(1),\dots,\gamma(l)$ with 
$\gamma(i)$ minimally adjacent to $\gamma(i+1)$ relating two such minimal
representatives. By minimality of $\gamma(0)$ and $\gamma(l)$
there is a first index $i$ such that $\gamma(i+1)$ is an infinitesimal
contraction of $\gamma(i)$ along $[R,S]$.
This implies that the point $R$ or $S$ has been added 
before, thus contradicting the minimality of the sequence
$\gamma(0),\dots,\gamma(l)$. 
\end{proof}

\begin{rem}
\begin{enumerate}
\item Infinitesimal contractions admit a sort of inverse: 
Given a generic configuration
$\mathcal P$ of $n$ points, choose a point $R\in \mathcal P$ and one of the
$(n-1)$ opposite pairs of cones delimited by lines through $R$ and 
$\mathcal P\setminus\{R\}$. Add now a new point
$S$ very close to $R$ in the interior 
of one of the two chosen opposite cones. This clearly yields
a configuration with $n(R,S)=0$ and every 
configuration of $(n+1)$ points infinitesimally adjacent to $\mathcal P$
can be constructed in this way. 
\item The notion of infinitesimal equivalence behaves well
with respect to isotopy: If two generic configurations are
isomorphic but not isotopic,
then their minimal representatives are also only isomorphic but
not isotopic.
\end{enumerate}
\end{rem}

\section{Monochromatic configurations}\label{mono_config}
This section is devoted to the study of {\it monochromatic configurations}
which are configurations with trivial Orchard partition.  

\begin{figure}[h]
\epsfysize=5cm
\centerline{\epsfbox{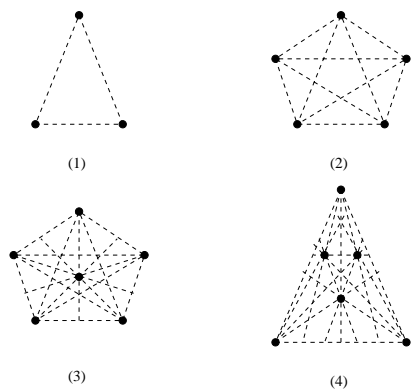}}
\caption{Monochromatic configurations with up to $6$ points}\label{up6points_mono}
\end{figure}

\begin{figure}[h]
\epsfysize=8cm
\centerline{\epsfbox{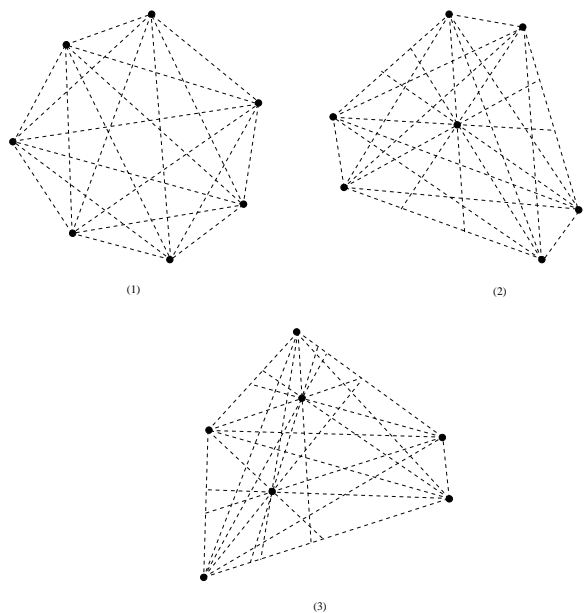}}
\caption{Monochromatic configurations of $7$ points having convex hull
  of size $5,6$ or $7$}\label{7points_mono-ch567}
\end{figure}

\begin{figure}[h]
\epsfysize=8cm
\centerline{\epsfbox{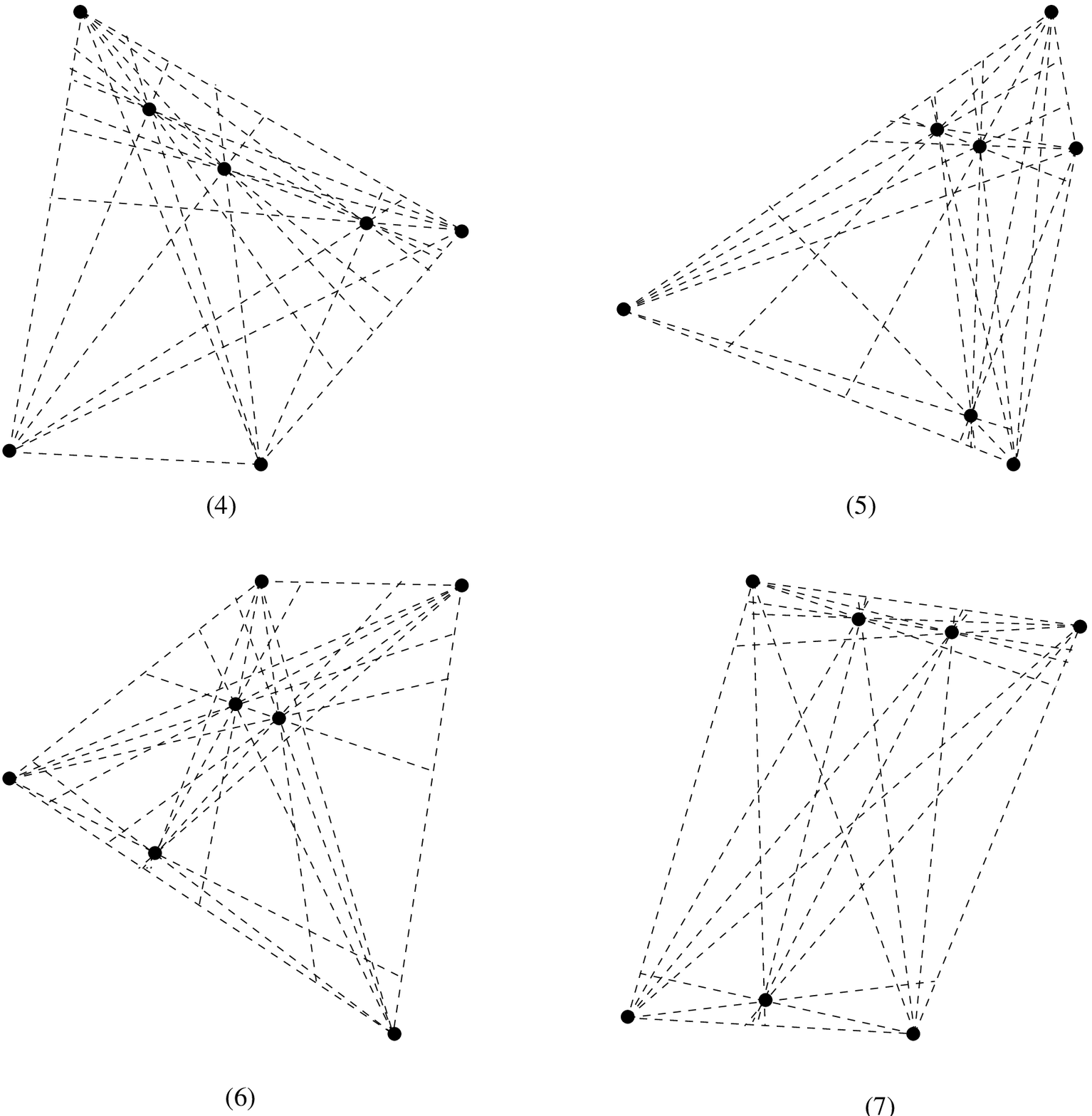}}
\caption{Monochromatic configurations of $7$ points with convex hull
  of size $4$}\label{7points_mono-ch4}
\end{figure}

\begin{figure}[h]
\epsfysize=8cm
\centerline{\epsfbox{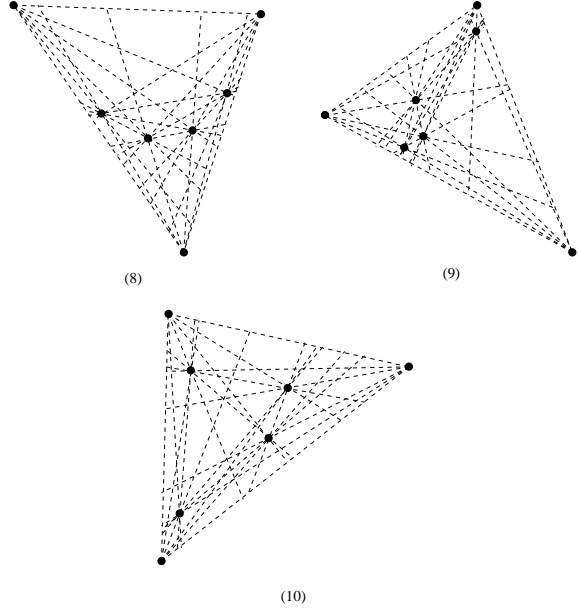}}
\caption{Monochromatic configurations of $7$ points with convex hull
  of size $3$}\label{7points_mono-ch3}
\end{figure}

Figures \ref{up6points_mono}-\ref{7points_mono-ch3} 
contain all monochromatic configurations (up
to non-oriented isomorphism) having at most $7$ points.

\newpage

\newpage

\section{Monochromatic constructions}\label{mono_construct}

This section contains a few constructions of
monochromatic configurations.
 
\begin{prop}
Let $\mathcal P$ be a monochromatic configuration having two 
infinitesimally-close points $P$ and $Q$ in its convex hull. 
Then $\mathcal P \setminus \{ P,Q \}$ is also monochromatic. 
\end{prop}

\begin{proof}
Use Proposition \ref{infi_close_prop}.
\end{proof}

This proposition can be used as follows: 

\begin{const}\label{infi_close_convex}
Let $\mathcal P=\mathcal P '\cup\{P,Q\}$ be a generic configuration
with an odd number of points having two infinitesimally-close vertices
$P,Q$ in its convex hull.
Suppose that $\mathcal P'$ is monochromatic and $n(P,P_0)$ is even for
some point $P_0\in \mathcal P '$. Then $\mathcal P$ is
monochromatic.
\end{const}

\begin{proof}
Apply the previous proposition.
\end{proof}

\begin{rem}
The same construction works also for a configuration $\mathcal P$ 
with an Orchard partition $(n,n)$: 
Add (generically) two infinitesimally-close points $P,Q$ 
to $\mathcal P$ such that $P,Q$ are in the convex hull 
of $\mathcal P \cup \{ P,Q \}$. The configuration
$\mathcal P \cup \{ P,Q \}$ 
has then Orchard partition $(n+1,n+1)$.
\end{rem}

There are of course many more constructions. Let us close this
section by mentioning another one:

Given a monochromatic configuration $\mathcal P$ having 
an odd number of points,
add generically two infinitesimally-close points $P,Q$. This
yields a configuration $\mathcal P'$ having two classes
$\mathcal C ^1$ and $\mathcal C ^2$, with $P,Q\in\mathcal C^1$,
which are separated by the line $L$ through $P,Q$. Add now another
pair of
infinitesimally-close points $P',Q'$ defining a line $L'$ parallel
to $L$ in such a way that $L'$ separates $P,Q$ from all points in $
\mathcal C ^2$. In the resulting configuration, the classes
$\mathcal C ^1$ and $\mathcal C^2$ merge. Only the new points 
$P',Q'$ might be outside this class. 
In this case, shift them along $L'$ into a suitable position.

\section{Some families of monochromatic configurations}\label{mono_families}

In this section we describe some monochromatic families.

\begin{prop}\label{convex_hull_family}
Vertices of a convex polygon having $2n+1$ edges yield a
monochromatic configuration. 
\end{prop}

\begin{proof}
Two adjacent vertices $P$ and $Q$ of the polygon
satisfy clearly $n(P,Q)=0$. This implies $P\sim Q$ since $n$ is odd.
\end{proof}

\begin{rem}
The above family can also be obtained by iterating Construction
\ref{infi_close_convex} starting with the unique configuration
consisting of one point.
\end{rem}

\begin{rem}\label{even_polygon}
The same argument shows that vertices of a convex polygon with
$2n$ edges have alternating colours.
\end{rem}

Another monochromatic family can be described as follows:

\begin{prop}
The following configuration $\mathcal P$ of $2n$ points is
mono\-chromatic for $n\geq 3$: Consider $n$ vertices spanning a
regular polygon $R$ (any convex polygon will in fact work) having $n$ sides. 
Add  $n$ more points which are interior to $R$ and very close to
the midpoints of all edges in $R$
(Figure \ref{family2_example} displays the example $n=5$).

\begin{figure}[h]
\epsfysize=4cm
\centerline{\epsfbox{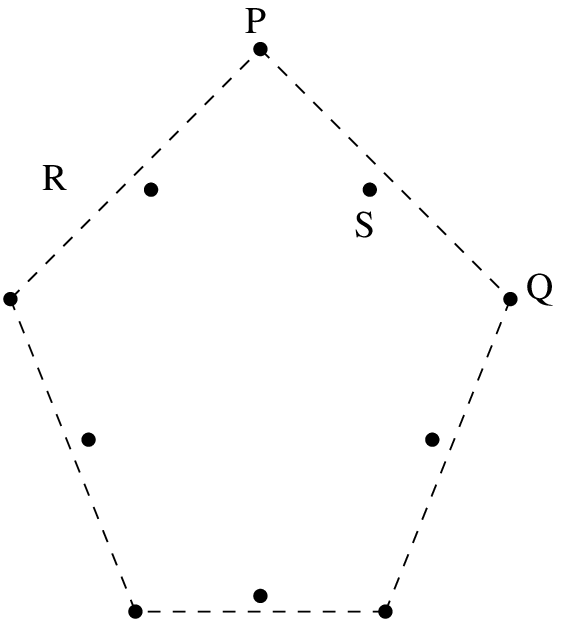}}
\caption{An example for a member in the family for $n=5$} \label{family2_example}
\end{figure}                                                                 

\end{prop} 

\begin{proof}
Consider first two adjacent vertices $P$ and $Q$ of the polygon $R$.
Denote by $S\in \mathcal P$ the point near the barycenter $(P+Q)/2$
of $P$ and $Q$. The points $P$ and $Q$ are separated by all lines
through $S$ and one of the remaining $2n-3$ points of ${\mathcal P}
\setminus\{P,Q,S\}$. Moreover, there are two more separating lines
going through vertices adjacent to $P$ and $Q$ and through
the corresponding near-midpoints of the edges adjacent to the
edge $[P,Q]$. We have hence $n(P,Q)=2n-1\equiv 2n-3\pmod 2$
which shows that $P\sim Q$.  

The points $P$ and $S$ are separated by exactly one of the two extra
lines mentioned above and hence we have $P\sim S$ too.
\end{proof}

\begin{rem}
One can prove the last result also by starting with a
convex polygon having $2n$ vertices $\{s_0,\dots,s_{2n-1}\}$.
By Remark \ref{even_polygon}, we have $s_a\equiv s_b$ if and
only if $a\equiv b\pmod 2$. Perform now flips with respect to all
triangles of vertices $s_{2i},s_{2i+1},s_{2i+2}$, $0 \leq i \leq n-1$ 
(with indices read $\pmod {2n}$). 
A vertex $s_k$ is hence involved in either $2$ or
$1$ flips according to the parity of $k$ and the 
resulting final configuration is monochromatic by the Flip
Proposition (Proposition \ref{flip_prop}).
\end{rem}

A monochromatic family involving congruences is given as follows:

\begin{prop}\label{family4}
Let $n \geq 5$ be an odd number, such that $(n+1) \equiv 0 \pmod 8$ or 
$(n+1) \equiv 6 \pmod 8$. The $n$ vertices of a regular polygon
having $n$ sides together with its barycenter yield
a monochromatic configuration $\mathcal P$ of $(n+1)$ points
(Figure \ref{family4_example} shows the case $n=5$).

\begin{figure}[h]
\epsfysize=4cm
\centerline{\epsfbox{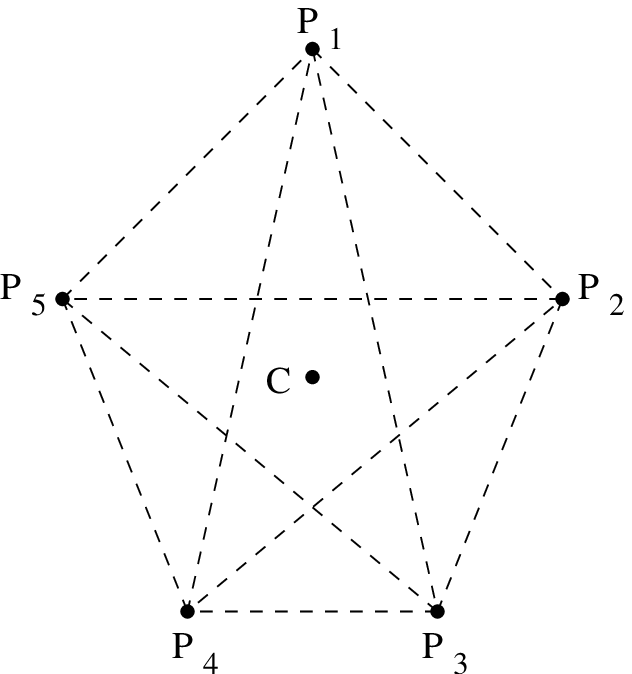}}
\caption{An example for a member in the family for $n=5$} \label{family4_example}
\end{figure}                                                                 

\end{prop} 

\begin{proof}
It is easy to check that 
$n(P,Q)=1$ for two adjacent vertices of the regular $n-$gon. 

We have yet to show that $n(P,C)$ is odd when $P$ is a vertex 
of the regular $n-$gon and $C$ is its barycenter. Any line $L$ separating
$C$ from $P$ comes together with its mirror obtained by reflecting
$L$ with respect to the line joining $C$ and $P$
(see lines $L_1$ and $L_2$ in Figure \ref{family4_proof_fig1}). 
Only separating lines orthogonal to the segment $[C,P]$ 
(like for instance line $L_3$ in Figure \ref{family4_proof_fig1})
have hence to be taken 
into acount when counting separating lines $\pmod {2}$. 

\begin{figure}[h]
\epsfysize=5cm
\centerline{\epsfbox{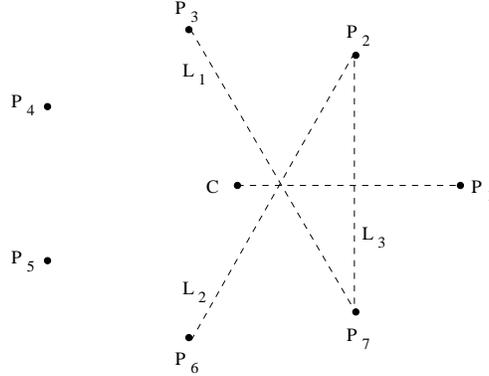}}
\caption{Examples for some separating lines for the 
  segment $[P_1,C]$, where $n=7$} \label{family4_proof_fig1}
\end{figure}                                                                 

Such lines are in bijection with the set
$$\{k \in \N \  \vert\
0<k\frac{2\pi}{n}<\frac{\pi}{2}\}=\{1,2,\dots,\lfloor
\frac{n}{4}\rfloor\}\ ,$$ 
and the result follows.
\end{proof}

\begin{rem}
For $n \equiv 1 \pmod 8$ or $n \equiv 3 \pmod 8$, the configurations considered in Proposition \ref{family4} have Orchard partition 
$(n,1)$.
\end{rem}

A somewhat similar family is given by:

\begin{prop}\label{family3}
Let $n \geq 3$ be an odd number. 
The following construction yields a monochromatic configuration of 
$2n+1$ points: Draw the $n$ vertices $\tilde P_1, \cdots, \tilde P_n$ 
of a regular $n$-gon.
Now replace each vertex $\tilde P_i$ by two infinitesimally-close
points $P_i,P'_i$, 
in order to obtain a configuration with convex hull of size $2n$. 
Add now a last point $C$ at the barycenter of this $2n$-gon
(See Figure \ref{7points_mono-ch567} number (2) for an example for $n=3$).
\end{prop} 

\begin{proof}
We have to check that $n(P_i,P'_i)$, 
$n(P'_i,P_{i+1})$ and $n(P_1,C)$ are even.

It is obvious  that $n(P_i,P'_i)=0$. Moreover,
$n(P'_i,P_{i+1})=2$ with two separating
lines defined by
the midpoint $C$ and the two vertices opposite to the
segment $[P_i',P_{i+1}]$.
 
It remains to check that $n(C,P_i)$ is even. This follows from the
fact that each line separating $C$ from $P_i$ has at least one of its 
defining points in a pair of infinitesimally-close points different
from $P_i$ and $P_i'$.
\end{proof}

\section{Some statistics}\label{statistics}

Aichholzer, Aurenhammer and Krasser have computed a data-base of all
generic configurations having $n\leq 10$ points, up to
unoriented isomorphism (see \cite{aak} and \cite{ak}). 
Assuming completeness of their data-base, we have computed the
Orchard partitions for all generic configurations having $n\leq 9$ points. 
The results for configurations of $7,8$ and $9$ points
are presented in Tables \ref{table1}-\ref{table3} according to the size of the 
configuration's convex hull and the size of the two
equivalence classes of the Orchard partition. 

Monochromatic configurations are of course partitions of type $(0,*)$.

\medskip

\begin{center}

\begin{table}[h]
\begin{tabular}{|c|r|r|r|r|}
\hline
Partition:          & (0,7)  & (1,6)  &  (2,5) &  (3,4) \\   
                    & mono   &        &        &        \\
Convex hull         &        &        &        &        \\
\hline
   3                &    3   &    7   &   13   &   26 \\
   4                &    4   &   11   &   16   &   28 \\
   5                &    1   &    3   &    5   &   13 \\
   6                &    1   &    0   &    1   &    2 \\
   7                &    1   &    0   &    0   &    0 \\
\hline
{\rm Sum:}          &   10   &   21   &   35   &   69 \\ 
\hline
\end{tabular}
\caption{Statistics for configurations with $7$ points}\label{table1}
\end{table}

\medskip

\begin{table}[h]
\begin{tabular}{|c|r|r|r|r|r|}
\hline
Partition:         & (0,8)   & (1,7) &  (2,6) &  (3,5) &  (4,4) \\   
                   & mono    &       &        &        &         \\
Convex hull        &         &       &        &        &          \\
\hline
   3               &   10    &   38  &   252  &   552  &   326 \\
   4               &   12    &   92  &   323  &   635  &   406 \\
   5               &    4    &   29  &    87  &   261  &   189 \\
   6               &    1    &    4  &    11  &    38  &    36 \\
   7               &    1    &    1  &     0  &     3  &     3 \\
   8               &    0    &    0  &     0  &     0  &     1 \\
\hline
{\rm Sum:}         &   28    &  164  &   673  &  1489  &   961 \\
\hline
\end{tabular}
\caption{Statistics for configurations with $8$ points}\label{table2}
\end{table}

\medskip

\begin{table}[h]
\begin{tabular}{|c|r|r|r|r|r|}
\hline
Partition:         & (0,9)   &  (1,8) &  (2,7) &  (3,6) &  (4,5)  \\   
                   & mono    &        &        &        &         \\
Convex hull        &         &        &        &        &         \\
\hline
   3               &  272    &  2469  &  8459  & 17493  & 26542 \\
   4               &  306    &  2484  & 10012  & 23234  & 34439 \\
   5               &  231    &  1277  &  4184  &  9273  & 13267 \\
   6               &   52    &   230  &   661  &  1490  &  2119 \\
   7               &    6    &    13  &    42  &   102  &   148 \\
   8               &    1    &     0  &     3  &     3  &     4 \\
   9               &    1    &     0  &     0  &     0  &     0 \\
\hline
{\rm Sum:}         &  869    &  6473  & 23361  & 51595  & 76519 \\
\hline
\end{tabular}
\caption{Statistics for configurations with $9$ points}\label{table3}
\end{table}

\end{center}

\section*{Acknowledgment} 
The second author wishes to thank Mikhail
Zaidenberg  and the Institut Fourier for hosting his stay.

\end{document}